\documentclass[11pt,a4paper]{article}
\usepackage[utf8]{inputenc}
\usepackage[english]{babel}
\usepackage{graphicx}
\usepackage{subcaption}
\usepackage{float}
\usepackage{pdfpages}
\usepackage{subcaption}
\usepackage{amsthm}
\usepackage{url}
\usepackage{xcolor}
\newtheorem{theorem}{Theorem}[section]

\newtheorem{assumption}{Assumption}[section]

\theoremstyle{definition}
\newtheorem{exmp}{Example}[section]

\theoremstyle{definition}

\usepackage{authblk}
\usepackage{amsmath}
\usepackage{amsfonts}
\usepackage{amssymb}
\usepackage{graphicx}
\usepackage{algorithmic}
\usepackage[linesnumbered,ruled,vlined]{algorithm2e}
\numberwithin{equation}{section}
\usepackage{lipsum}
\usepackage{mathtools}
\usepackage{scalefnt}
\usepackage{empheq}
\usepackage[a4paper,left=3.5cm,right=3cm,top=3cm,bottom=3cm]{geometry}
\usepackage{makeidx}
\usepackage{bm}
\makeindex

\usepackage{indentfirst}

\date{}
\title{An initial-corrected splitting approach for convection-diffusion-reaction problems}

\author[1]{Thi Tam Dang\thanks{tam.dang@helsinki.fi}}
\author[2]{Lukas Einkemmer\thanks{lukas.einkemmer@uibk.ac.at}}
\author[2]{Alexander Ostermann\thanks{alexander.ostermann@uibk.ac.at}}
\affil[1]{Department of Mathematics and Statistics, University of Helsinki, Finland}
\affil[2]{Department of Mathematics, University of Innsbruck, Austria}

\begin{document}

\maketitle

\begin{abstract}
Splitting methods constitute a well-established class of numerical schemes for solving convection-diffusion-reaction problems. They have been shown to be effective in solving problems with periodic boundary conditions. However, in the case of Dirichlet boundary conditions, order reduction has been observed even with homogeneous boundary conditions.
In this paper, we propose a novel splitting approach, the so-called \emph{initial-corrected splitting method}, which succeeds in overcoming order reduction. A convergence analysis is performed to demonstrate second-order convergence of this modified Strang splitting method. Furthermore, we conduct numerical experiments to illustrate the performance of the newly developed splitting approach.
\end{abstract}

\section{Introduction}

In this paper, we consider the following abstract convection-diffusion-reaction problem
\begin{equation}\label{B1.1}
	\begin{aligned}
		&	\partial_{t} u(t) =  D u(t) + a \cdot \nabla u(t) + f(t,u(t))  , \\
		& u(t)|_{\partial \Omega}= b(t),\\
		&u(0)= u_0,
	\end{aligned}
\end{equation}
where $D$ denotes a second-order elliptic differential operator and $a$ is a space-dependent convection coefficient. The nonlinear source term $f$ is assumed to be sufficiently smooth. We consider this problem on a bounded domain $\Omega \subset \mathbb{R}^{d}$ with smooth boundary $\partial \Omega$, initial data $u_0$, and Dirichlet boundary conditions given by $b: [0,T]\times \partial \Omega \rightarrow \mathbb{R}$.  The boundary data $b$ may be time-dependent and is assumed to be sufficiently smooth.

In the context of convection-diffusion-reaction problems, the choice of boundary conditions is a significant factor in determining the overall behavior of the solutions. For periodic boundary conditions, there are usually no issues, and no special care has to be taken in the analysis. However, Dirichlet conditions are often needed in applications. In this case, a specific value of the solution at the boundary is imposed.

The pursuit of efficient numerical methods for solving convection-diffusion-reaction problems is an active area of research. Among the various methods available, splitting methods are widely studied and can be used to efficiently solve such problems. In this approach we split the convection-diffusion-reaction problem into two separate subproblems:
\begin{equation}\label{B1.2}
	\begin{aligned}
		&\partial_{t} v(t) = D v(t) + f(t,v(t)), \\
		&	v(t)|_{\partial \Omega}= b(t),\\	
	\end{aligned}
\end{equation}	
and
\begin{equation}\label{B1.3}
	\begin{aligned}
		&	\partial_{t} w(t)= a \cdot  \nabla w(t),\\
		&	w(t)|_{\Gamma} = c(t),
	\end{aligned}
\end{equation}	
where $\Gamma$ denotes the inflow boundary and $c(t)=b(t)\vert_{\Gamma}$.

The primary advantage of splitting methods is that they allow for the independent treatment of each subproblems, i.e.~they do not require a monolithic solver for the full problem \eqref{B1.1}. Instead, specifically tailored methods for the subproblems can be used. For example, the diffusion-reaction problem \eqref{B1.2} can be solved efficiently by using implicit or IMEX methods and a multigrid preconditioner (among many other options). On the other hand, for the convection part an explicit scheme is usually sufficient. Furthermore, splitting methods can be very effective in preserving certain properties of the solution. For example, the numerical solution is positive if the numerical solution of the subproblems have this property (see, e.g., \cite{HANSEN20121428}).

Our study focuses on the Strang splitting method for solving the convection-diffusion-reaction problem \eqref{B1.1}.  Let $\tau$ be the time step size. Then, one step of the classical Strang splitting in the time interval $\left[ t_{n}, t_{n+1}\right]$ is given by
\begin{align}\label{B1.4}
	u_{n+1} = \varphi_{\frac{\tau}{2}}^{D,f} \circ \varphi_{\tau}^{w}\circ \varphi_{\frac{\tau}{2}}^{D,f}(u_{n}),
\end{align}
where $\varphi_{t}^{D,f}$ and $\varphi_{t}^{w}$ are the flows of the problems specified in \eqref{B1.2} and \eqref{B1.3}, respectively. For the purpose of our analysis, we assume that these flows are computed exactly. This corresponds to the situation where the splitting error dominates the numerical error for the subflows.

Strang splitting is a well-known second-order method (see, e.g., \cite{HUNDSDORFER1995191,mclachlan_quispel_2002}). However, numerical experiments have shown that even when homogeneous Dirichlet boundary conditions are imposed, the Strang splitting method of the form \eqref{B1.4} suffers from order reduction (see the experiments in Sections~4 and~5). Commonly only first order convergence is observed in this setting.

The reason for the order reduction is that when \eqref{B1.1} is split into \eqref{B1.2} and \eqref{B1.3}, Dirichlet boundary conditions are imposed on the diffusion-reaction equation \eqref{B1.2}, while inflow boundary conditions are imposed on the convection equation \eqref{B1.3}. As a result, during the internal step of the splitting process, there is a lack of compatibility between inflow boundary conditions and the prescribed Dirichlet boundary conditions. This observation encourages us to find modified splitting methods that do not suffer from order reduction.

Many researchers have devoted their efforts to adapting splitting methods to avoid the order reduction phenomenon that occurs in diffusion-reaction problems with nontrivial boundary conditions (see, e.g., \cite{10.1016/j.cam.2019.02.023,doi:10.1137/19M1257081, EINKEMMER201876, doi:10.1137/140994204,Einkemmer2016ACO, doi:10.1137/16M1056250,hundsdorfer2013numerical,r_1981}. These approaches, however, mainly focus on splitting between diffusion and reaction.  In this paper, we present a new splitting approach, the so-called \emph{initial-corrected splitting method}, which aims at overcoming the order reduction when splitting between diffusion, reaction and convection.

The main ideas of our proposed splitting method is to minimize the influence of the inflow boundary conditions by subtracting the initial data from the solution. This will effectively eliminate the incompatibility condition that is encountered in the internal step of the classical Strang splitting scheme. We then rewrite the convection-diffusion-reaction problem \eqref{B1.1} so that homogeneous Dirichlet boundary conditions are imposed and the initial data are zero. The detailed construction of our splitting schemes will be described in Section~2.

Our new splitting approach is constructed directly from the initial data, requiring no computation of a correction term (as is common in the literature; see, e.g., \cite{10.1016/j.cam.2019.02.023,doi:10.1137/19M1257081,doi:10.1137/140994204,doi:10.1137/16M1056250}). This leads to an accurate scheme that is also computationally efficient.

In addition to proposing the initial-corrected splitting method, we also conduct a convergence analysis of the proposed method applied to \eqref{B1.1}. We prove that the method based on Strang splitting achieves second-order convergence. This analysis provides further insight into the convergence properties of our method. Furthermore, a series of numerical experiments will be performed to illustrate the performance of our new splitting approach and to show its versatility for a wide range of problems.

The paper is organized as follows: Section~2 describes the construction of the initial-corrected splitting method. In Section~3, we provide a comprehensive error analysis of the modified Strang splitting scheme, including a detailed proof of second-order convergence. Finally, Sections~4 and~5 present a series of numerical experiments in one and two dimensions.

\section{Description of the method}	

In this section, we will provide a detailed description of the newly proposed splitting approach, referred to as the initial-corrected splitting method, for the convection-diffusion-reaction problem \eqref{B1.1}. Let us first consider the case of time-independent boundary conditions. By doing so we aim to present a clear and coherent explanation of the main ideas underlying our approach.

Our goal is to rewrite problem \eqref{B1.1} in such a way that homogeneous Dirichlet boundary conditions are imposed and the initial data are set to zero. The transformation turns out to be solution-dependent. Let $u_n$ denote the numerical solution of \eqref{B1.1} at time $t=t_n$ and $\tau$ the (actual) step size of the integrator. In order to determine the numerical approximation $u_{n+1}$ at time $t_{n+1} = t_n+\tau$, we have to solve \eqref{B1.1} with initial data $u(t_n) = u_n$. For this we set
\begin{equation}\label{eq:z-invariant}
z_n = u_n
\end{equation}
and $\hat{u}(t) = u(t)-z_n$ for $t_n\le t \le t_{n+1}$. Then $\hat{u}(t)$ satisfies the following initial boundary value problem
\begin{subequations}\label{B2.4}
	\begin{align}
		&\partial_{t} \hat{u}(t) = D\hat{u}(t) + a \cdot \nabla \hat{u}(t) + h(t,\hat{u}(t)), \label{eq:2.4a} \\
		& \hat{u}|_{\partial \Omega}= 0, \label{eq:2.4b} \\
		&\hat{u}(t_n)= \hat{u}_n =0, \label{eq:2.4c}
	\end{align}
\end{subequations}
where $h$ denotes the following modified nonlinearity
\begin{align}\label{B2.5}
	h(t,\hat{u}(t))= f(t, \hat{u}(t) + z_n) + D z_n + a\cdot\nabla z_n.
\end{align}
By setting $\hat{v}(t)= v(t)-z_n $,\: $\hat{w}(t)= w(t)-z_n $, we now split \eqref{B2.4} as follows: 
\begin{equation}\label{B2.6}
	\begin{aligned}
		&	\partial_{t} \hat{v}(t)=  D \hat{v}(t) + h(t,\hat{v}(t)),\\
		&	\hat{v}|_{\partial \Omega}= 0,
	\end{aligned}
\end{equation}
and
\begin{equation}\label{B2.7}
	\begin{aligned}
		&	\partial_{t} \hat{w}(t)= a \cdot \nabla \hat{w}(t) ,\\
		& \hat{w}|_{\Gamma} = 0.
	\end{aligned}
\end{equation}
We call this splitting a modified (Strang) splitting or more precisely an initial-corrected splitting scheme. One step of this initial-corrected Strang splitting scheme for time-invariant boundary conditions is described in Algorithm \ref{alg 1} below.

\medskip
\begin{algorithm}[H]\label{Al1}
	\SetAlgoLined
	\KwIn{Initial data $u_{n}$ at time $t_n$, step size $\tau$}
	\KwOut{Numerical solution $u_{n+1}$ at time $t_{n+1}$}
    Let $z_n = u_n$ and $\hat{v}(t_n) = u_n-z_n = 0$\;
	Compute the solution of \eqref{B2.6} with initial data $\hat{v}(t_n)$ to obtain $\hat{v}\left(t_n+\frac{\tau}{2} \right)$\;
	Compute the solution of \eqref{B2.7} with initial data $\hat{w}(t_n)= \hat{v}\left(t_n+\frac{\tau}{2} \right)$ to obtain $\hat{w}(t_{n}+\tau)$\;
	Compute the solution of \eqref{B2.6} with initial data $\hat{v}\left( t_n+\frac{\tau}{2}\right)=  \hat{w}(t_n+\tau)$ to obtain $\hat{u}_{n+1}= \hat{v}(t_{n}+\tau)$\; 
	Compute $u_{n+1}= \hat{u}_{n+1} + z_n$\;
	\caption{Initial-corrected Strang splitting for \eqref{B1.1}; time-invariant $b$}
	\label{alg 1}
\end{algorithm}

\medskip

For time-dependent boundary conditions, i.e., $b = b(t)$, additional challenges arise. The piecewise constant correction $z_n = u_n$ is no longer accurate enough and it is necessary to provide a time-dependent correction $z_n(t)$. Since the boundary values of $z_n(t)$ should approximate $b(t) = u(t)|_{\partial \Omega}$, the simplest choice for the correction is a first-order Taylor approximation. Therefore, we set
\begin{align}\label{2.7}
	z_n(t) = u_{n} + (t-t_{n}) \big( Du_{n}+ a \cdot \nabla u_{n}+ f(t_{n}, u_{n})\big)
\end{align}
and carry out the transformation $\hat{u}(t) = u(t)-z_n(t)$ for $t_n\le t \le t_{n+1}$. In this way, get again \eqref{B2.4}, however, the nonlinearity $h$ now has the form
\begin{equation}\label{2.8}
\begin{aligned}
	  h(t, \hat{u}(t)) &= f(t, \hat{u}(t) + z_n(t)) - f(t_n, u_n)\\
 &\quad + (t-t_n)\bigl(D + a\cdot \nabla\bigr)\bigl(D u_n + a\cdot\nabla u_n + f(t_n, u_n)\bigr).
\end{aligned}
\end{equation}
This leads to the following algorithm for time-dependent boundary conditions.

\medskip
\begin{algorithm}[H]\label{Al2}
	\SetAlgoLined
	\KwIn{Initial data $u_{n}$ at time $t_n$, step size $\tau$}
	\KwOut{Numerical solution $u_{n+1}$ at time $t_{n+1}$}
    Let $z_n(t)$ and $h(t,\hat u(t))$ be given by \eqref{2.7} and \eqref{2.8}, respectively, and let $\hat{v}(t_n) = u_n-z_n(t_n) = 0$\;
	Compute the solution of \eqref{B2.6} with initial data $\hat{v}(t_n)$ to obtain $\hat{v}\left(t_n+\frac{\tau}{2} \right)$\;
	Compute the solution of \eqref{B2.7} with initial data $\hat{w}(t_n)= \hat{v}\left(t_n+\frac{\tau}{2} \right)$ to obtain $\hat{w}(t_{n}+\tau)$\;
	Compute the solution of \eqref{B2.6} with initial data $\hat{v}\left( t_n+\frac{\tau}{2}\right)=  \hat{w}(t_n+\tau)$ to obtain $\hat{u}_{n+1}= \hat{v}(t_{n}+\tau)$\; 
	Compute $u_{n+1}= \hat{u}_{n+1} + z_n(t_n+\tau)$\;
	\caption{Initial-corrected Strang splitting for \eqref{B1.1}; time-dependent $b$}
	\label{alg 2}
\end{algorithm}

\medskip

Our proposed splitting approach offers a notable advantage over the methods presented in \cite{10.1016/j.cam.2019.02.023,doi:10.1137/19M1257081,doi:10.1137/140994204,doi:10.1137/16M1056250}  when dealing with time-dependent boundary conditions. Our approach eliminates the need to compute a correction term, and $z_n$ or $z_n(t)$ can be computed directly from the initial data for each step in advance, greatly reducing both implementation complexity and computational cost.
In the next section a comprehensive convergence analysis is performed to theoretically establish the second-order convergence property of the proposed splitting method.

\section{Convergence analysis}
This section is devoted to conduct a thorough error analysis of the modified Strang splitting method applied to \eqref{B1.1}. First, we provide the analytical framework that we use extensively in our convergence analysis. Then, we will study both the local error and the global error of the modified Strang splitting method applied to \eqref{B1.1}.
\subsection{Analytical framework}
In this section, we will consider the following abstract evolution equation
\begin{equation}\label{B3.1}
	\begin{aligned}
		&\partial_{t}\hat{u}(t) +A \hat{u}(t)= a  \cdot \nabla \hat{u}(t)+ h(t,\hat{u}(t)),
	\end{aligned}
\end{equation}
where $h(t,\hat{u}(t))$ is defined as in \eqref{B2.5} or \eqref{2.8}. The linear operator $A$ is defined as $Av = -D v$ for all $v \in D(A)$. For instance, if we consider the problem in $L^2(\Omega)$, the domain $D(A)$ of $A$ is $H^2(\Omega) \cap H_0^1(\Omega)$, for $D$ being a second-order strongly elliptic differential operator. Note that the problem \eqref{B3.1} is equivalent to the problem \eqref{eq:2.4a}--\eqref{eq:2.4b}, with the homogeneous Dirichlet boundary conditions incorporated into the domain of the operator $A$. For the purpose of error analysis, we consider the following assumptions:
\begin{assumption}\label{Bas1}
	Let $X$ be a Banach space  equipped with the norm $\Vert \cdot \Vert$. We assume that the operator $ -A$ is the infinitesimal generator of an analytic semigroup $ e^{-tA}$. 	
\end{assumption}
There exists a constant $\omega \geq 0$ such that fractional powers $(\omega I +A)^{\gamma}$ are well-defined for $\gamma \in \mathbb{R}$, see \cite{henry1981geometric}. Let us recall some properties of the analytic semigroup $e^{-t A}$ that  will be extensively used throughout this section. For $t\in [0,T]$ there exists a constant $C$ such that
\begin{align}\label{B3.2}
	\Vert e^{-tA}\Vert \le C.
\end{align}
Moreover, the analytic semigroup $e^{-t A}$ satisfies the following parabolic smoothing property
\begin{align}\label{B3.3}
	\| (\omega I+A)^{\gamma} e^{-tA}\| \le C t^{-\gamma},\qquad\gamma>0,
\end{align}
uniformly for $t \in (0,T]$.

We need an appropriate framework to incorporate the nonlinearity $f$, which we will establish in the next assumption:
\begin{assumption}\label{Bas2}
	Let $f: X\rightarrow X $ be sufficiently Fr\'echet differentiable in a strip along the exact solution. We assume that all derivatives of $f$ are uniformly bounded.
\end{assumption}
In particular, Assumption \ref{Bas2} implies that $f$ is locally Lipschitz continuous in a strip along the exact solution $u(t)$. Consequently, there exists a constant $L(R)$ such that
\begin{align}\label{B3.4}
	\Vert f(t,u(t))- f(t, v)\Vert \le L \Vert u(t)-v\Vert,
\end{align}
for $\max \Vert u(t)-v\Vert \le R$.

The exact solution of \eqref{B3.1} at the time $t_{n+1} = t_{n}+\tau$ is given by the variation-of-constants formula as follows:
\begin{equation}\label{B3.5}
	\begin{aligned}
		\hat{u}(t_{n+1})= e^{-\tau A} \hat{u}(t_{n})+ \int_0^{\tau} e^{-(\tau-s)A} \Big( a \cdot \nabla \hat{u}(t_{n}+s)  + h(t_{n}+s, \hat{u}(t_{n}+s)) \Big) ds.
	\end{aligned}
\end{equation}
Therefore, the exact solution of \eqref{B1.1} can be expressed in the following way
\begin{equation}\label{B3.6}
	\begin{aligned}
		u(t_{n+1})&= z(t_{n+1}) +e^{-\tau A} \hat{u}(t_{n})+  \int_{0}^{\tau} e^{-(\tau-s)A} \Big( a \cdot \nabla \hat{u}(t_{n}+s) \Big)ds \\
		& \qquad +  \int_{0}^{\tau} e^{-(\tau-s)A} h(t_{n}+s, \hat{u}(t_{n}+s)) ds,
	\end{aligned}
\end{equation}
where $z(t) = u(t)-\hat{u}(t)$ is the correction employed in the transformation. 

We now proceed by splitting \eqref{B3.1} into the following two subproblems
\begin{equation}\label{B3.7}
	\begin{aligned}
		\partial_{t} \hat{v}(t) + A v(t)= h(t, \hat{v}(t)),	
	\end{aligned}
\end{equation}
and
\begin{equation}\label{B3.8}
	\begin{aligned}
		&	\partial_{t} \hat{w}(t) = a \cdot \nabla \hat{w}(t),\\
		&  \hat{w}(t)|_{\Gamma}= 0.
	\end{aligned}
\end{equation}

Our main convergence result for the modified Strang splitting applied to \eqref{B1.1} is stated in Theorem~\ref{Th: 2.1}. For that, we employ the following assumption on the data of~\eqref{B1.1}.
\begin{assumption}\label{Bas3}
	Let $\Omega$ be a bounded domain in $\mathbb{R}^{d}$, let $D$ be a second-order strongly elliptic differential operator with smooth coefficients, and let the space-dependent convection coefficient $a$ and the boundary data $b(t)$ be sufficiently smooth. We further assume that the solution $u(t)$ of \eqref{B1.1} is sufficiently smooth.
\end{assumption}

We are now in a position to formulate our main convergence result.
\begin{theorem} \label{Th: 2.1}
	Let the Assumptions \ref{Bas1}-\ref{Bas3} be satisfied. Then there exists a constant $\tau_0>0$ such that for all step sizes $0<\tau \le \tau_0$ and $t_{n}=n\tau$ we have that the modified Strang splitting applied to \eqref{B1.1} satisfies the global error bound
	\begin{align}
		\Vert u_{n}-u(t_{n}) \Vert \le C\tau^2(1 + \left| \log \tau \right| ), \ \ 0\le n\tau \le T,
	\end{align}
	where the constant $C$ depends on $T$ but is independent of $\tau$ and $n$.
\end{theorem}
In the following subsections, we will study the local error and the global error of the modified Strang splitting method applied to \eqref{B1.1} within the framework of analytic semigroups.

\subsection{Local error}
In this subsection, we derive the local error bound for the modified Strang splitting applied to \eqref{B1.1}. For this purpose, let us consider one step of the numerical solution, starting at time $t_{n}$ with the initial value $u(t_n)$ on the exact solution. The solution of \eqref{B3.7} with the initial value $\hat{v}(t_{n})= \hat{u}(t_{n}) = u(t_n) - u_n$ and step size $\tau/2$ can be expressed by the variation-of-constants formula
\begin{equation}\label{B3.10}
	\begin{aligned}
		Y= \hat{v}\left( t_{n}+ \frac{\tau}{2}\right)= e^{-\frac{\tau}{2} A} \hat{u}(t_{n})+	\int_0^{\frac{\tau}{2}} e^{ -\left( \frac{\tau}{2}-s \right)A } h(t_{n}+s, \hat{v}(t_{n}+s)) ds.
	\end{aligned}
\end{equation}
The solution of \eqref{B3.8} with the full step size is represented by a Taylor series as follows
\begin{equation}\label{B3.11}
	\begin{aligned}
		\hat{w}(t_{n}+ \tau)  =  Y + \tau a \cdot \nabla Y + \frac{\tau^2}{2} a \cdot \nabla \left(  a\cdot \nabla Y \right) + \mathcal{O}(\tau^3),
	\end{aligned}
\end{equation}
where $\mathcal{O}(\tau^3)$ denotes the bounded remainder term of order $3$. Carrying out again a half step of \eqref{B3.7} with the initial value $\bar{v}\left( t_{n}+\frac{\tau}{2}\right)= \hat{w}(t_{n}+ \tau)$ yields one step of the modified Strang splitting scheme applied to \eqref{B1.1}:
\begin{equation}\label{B3.14}
	\begin{aligned}
		S_{\tau}u(t_{n})& = z(t_{n+1}) +e^{-\tau A} \hat{u}(t_{n})+  \int_0^{\frac{\tau}{2}} e^{-(\tau-s)A} h(t_{n}+s, \hat{v}(t_{n}+s)) ds \\
		& \qquad + \tau e^{-\frac{\tau}{2} A} a \cdot \nabla Y + \frac{\tau^2}{2} e^{-\frac{\tau}{2} A} a \cdot \nabla \left( a \cdot \nabla Y \right)  \\
		& \qquad + \int_{\frac{\tau}{2}}^{\tau} e^{-(\tau-s)A} h(t_{n}+s, \bar{v}(t_{n}+s)) ds + \mathcal{O}(\tau^3),
	\end{aligned}
\end{equation}
where again $z(t) = u(t)-\hat{u}(t)$ is the correction employed in the transformation. 

Let us denote by $\delta_{n+1}= \mathcal{S}_{\tau}u(t_{n})- u(t_{n+1})$ the local error. Subtracting the exact solution \eqref{B3.6} from the numerical solution \eqref{B3.14} gives the following representation of the local error
\begin{align}\label{B3.15}
	\delta_{n+1}=\delta_{n+1}^{\left[1 \right] } + \delta_{n+1}^{\left[2 \right] },
\end{align}
where
\begin{equation}\label{B3.16}
	\begin{aligned}
		&\delta_{n+1}^{\left[1 \right] } = \tau e^{-\frac{\tau}{2} A} a \cdot \nabla Y + \frac{\tau^2}{2} e^{-\frac{\tau}{2} A} a \cdot  \nabla \left( a \cdot \nabla Y \right) -\int_0^{\tau} e^{-(\tau-s)A} a \cdot \nabla \hat{u}(t_{n}+s) ds,  \\
	\end{aligned}
\end{equation}
and
\begin{equation}\label{B3.17}
	\begin{aligned}
		\delta_{n+1}^{\left[2 \right] } &= \int_0^{\frac{\tau}{2}} e^{-(\tau-s)A} h(t_{n}+s, \hat{v}(t_{n}+s)) ds \\
		& \qquad+ \int_{\frac{\tau}{2}}^{\tau} e^{-(\tau-s)A} h(t_{n}+s, \bar{v}(t_{n}+s)) ds\\
		& \qquad - \int_0^{\tau} e^{-(\tau-s)A} h(t_{n}+s, \hat{u}(t_{n}+s)) ds .
	\end{aligned}
\end{equation}


Let us first consider $\delta_{n+1}^{\left[1 \right]}$. By setting
\begin{align}\label{B3.18}
	k_{n}(s) = e^{-(\tau-s)A} \hat{k}_{n}(s), \quad  \hat{k}_{n}(s) =  a \cdot \nabla \hat{u}(t_{n}+s),
\end{align}
and using the midpoint rule, we get
\begin{equation}\label{B3.19}
	\begin{aligned}
		&\int_0^{\tau} k_{n}(s)ds = \tau k_{n}\left( \frac{\tau}{2}\right) +\int_0^{\tau}M(s, \tau) k_{n}^{\prime \prime}(s)ds,
	\end{aligned}
\end{equation}
with $M(s,\tau) = s^2/2$ for $s \le \tau/2$ and $M(s,\tau)= (\tau-s)^2/2$ for $s >\tau/2$. In view of \eqref{B3.18} and \eqref{B3.19}, we can rewrite $\delta_{n+1}^{\left[1\right]}$ as follows:
\begin{equation}\label{B3.20}
	\begin{aligned}
		\delta_{n+1}^{\left[1\right]}= \delta_{n+1}^{\left[1,1\right]} +\delta_{n+1}^{\left[1,2\right]},
	\end{aligned}
\end{equation}
where
\begin{align}
		&	\delta_{n+1}^{\left[1,1\right]} = \tau e^{-\frac{\tau}{2}  A} a \cdot  \bigl(  \nabla Y- \nabla \hat{u}\left(t_{n}+ \tfrac{\tau}{2} \right) \bigr) + \frac{\tau^2}{2} e^{-\frac{\tau}{2} A} a \cdot \nabla \left(  a \cdot  \nabla Y \right), \label{3.19}\\
	& \delta_{n+1}^{\left[1,2\right]} = - \int_0^{\tau}M(s, \tau) k_{n}^{\prime \prime}(s)ds. \label{3.20}
\end{align}
Using \eqref{B3.5} and \eqref{B3.10} we have
\begin{equation}\label{B3.22}
	\begin{aligned}
	\hat{u}\left( t_{n}+ \frac{\tau}{2} \right) - Y = \int_0^{\frac{\tau}{2}} \tilde{k}_{n}(s) ds + \mathcal{O}(\tau^2),
	\end{aligned}
\end{equation}
where $	\tilde{k}_{n}(s)$ is given by \eqref{B3.18} with $\tau$ replaced by $\tau/2$. Note that in the above estimate, we used $h(t_{n}+s, \hat{u}(t_{n}+s))- h(t_{n}+s, \hat{v}(t_{n}+s)) = \mathcal{O}(\tau)$. Taylor expansion of $\tilde{k}_{n}$ yields the formula
\begin{equation}\label{B3.23}
	\begin{aligned}
		\int_0^{\frac{\tau}{2}} \tilde{k}_{n}(s) ds &= \frac{\tau}{2} \tilde{k}_{n} \left( \frac{\tau}{2} \right) + \int_0^{\frac{\tau}{2}}\int_{\frac{\tau}{2}}^{s} \tilde{k}_{n}^{\prime}(\xi)d\xi ds.
	\end{aligned}
\end{equation}	
Plugging \eqref{B3.23} into \eqref{B3.22} we get
\begin{equation}\label{B3.24}
	\begin{aligned}
		\hat{u}\left( t_{n}+ \frac{\tau}{2} \right) - Y &=  \frac{\tau}{2} a \cdot \nabla \hat{u} \left( t_{n}+ \frac{\tau}{2} \right) + \int_0^{\frac{\tau}{2}}\int_{\frac{\tau}{2}}^{s}\tilde{k}_{n}^{\prime}(\xi)d\xi ds.
	\end{aligned}
\end{equation}
Inserting \eqref{B3.24} into \eqref{3.19} we obtain
\begin{equation}\label{B3.25}
	\begin{aligned}
		\delta_{n+1}^{\left[1,1\right]}  = - \tau e^{-\frac{\tau}{2}A} A A^{-1}a \cdot  \nabla  \int_0^{\frac{\tau}{2}} \int_{\frac{\tau}{2}}^{s}  \tilde{k}_{n}^{\prime}(\xi) d\xi ds+ \mathcal{O}(\tau^3).
	\end{aligned}
\end{equation}
To estimate $\delta_{n+1}^{ [1,1]}$, we need to bound the integral $I_{n}= \int_0^{\frac{\tau}{2}} \int_{\frac{\tau}{2}}^{s}  \tilde{k}_{n}^{\prime}(\xi) d\xi ds.$ Since the compatibility condition at the boundary is satisfied, i.e. $\nabla \hat{u}(t_{n})|_{\partial \Omega} =0$, it follows that $\hat{k}_{n}(0) \in D(A)$. Therefore
\begin{equation}\label{B3.26}
	\begin{aligned}
		\| I_{n}\| &= \left\| \int_0^{\frac{\tau}{2}} \int_{\frac{\tau}{2}}^{s}  \tilde{k}_{n}^{\prime}(\xi) d\xi ds \right\| \\
		& \le \int_0^{\frac{\tau}{2}} \int_{\frac{\tau}{2}}^{s} \| e^{-\left(\frac{\tau}{2}-\xi \right)A } A \hat{k}_{n}(\xi) + e^{-\left(\frac{\tau}{2}-\xi \right)A } \hat{k}_{n}^{\prime}(\xi) \| d\xi ds\\
		& \le \int_0^{\frac{\tau}{2}} \int_{\frac{\tau}{2}}^{s} \| e^{-\left(\frac{\tau}{2}-\xi \right)A } \| \| A \hat{k}_{n}(0) \| d\xi ds\\
		& \qquad +  \|A e^{-\frac{\tau}{4}A} \| \int_0^{\frac{\tau}{2}} \int_{\frac{\tau}{2}}^{s} \|e^{-\left(\frac{\tau}{4}-\xi \right)A} \| \| \mathcal{O}(\xi)\| d\xi ds\\
		&  \qquad + \int_0^{\frac{\tau}{2}} \int_{\frac{\tau}{2}}^{s}  \| e^{-\left(\frac{\tau}{2}-\xi \right)A } \| \| \hat{k}_{n}^{\prime}(\xi) \| d\xi ds\\
		& \le C\tau^2.
	\end{aligned}
\end{equation}
Thanks to the sufficient smoothness of $a$, it follows that $\|A^{-1} a \cdot \nabla\|$ is bounded. Using \eqref{B3.2}, \eqref{B3.25}, and \eqref{B3.26}, we get
\begin{align}\label{B3.28}
	\delta_{n+1}^{\left[1,1\right]} =A \mathcal{O}(\tau^3)+ \mathcal{O}(\tau^3).
\end{align}

To bound $\| \delta_{n+1}^{\left[1,2\right]}\| $, it is sufficient to show the boundedness of $k_{n}^{\prime \prime}(s)$. The second time derivative of $k_{n}(s)$ is given by
\begin{equation}
	\begin{aligned}
		k_{n}^{\prime \prime}(s)&= A  e^{-(\tau-s)A} A \hat{k}_{n}(s)+ 2 A e^{-(\tau-s) A} \hat{k}_{n}^{\prime}(s)
		+ e^{-(\tau-s) A} \hat{k}_{n}^{\prime \prime}(s),
	\end{aligned}
\end{equation}
where $\hat{k}_{n}(s)$ given in \eqref{B3.18}.
Thus, we have
\begin{equation}\label{B3.28a}
	\begin{aligned}
		\delta_{n+1}^{\left[1,2\right]} &= - \int_0^{\tau} M(s,\tau) A e^{-(\tau-s)A} A \hat{k}_{n}(s) ds - 2 \int_0^{\tau} M(s,\tau) A e^{-(\tau-s) A} \hat{k}_{n}^{\prime}(s) ds\\
		& \qquad - \int_0^{\tau} M(s,\tau) e^{-(\tau-s) A} \hat{k}_{n}^{\prime \prime}(s) ds\\
		& = -P(s)-Q(s)-R(s).
	\end{aligned}
\end{equation}
Let us now estimate the various integrals on the right-hand side of \eqref{B3.28a}.
We set
\begin{align*}
	P(s)&= \int_0^{\tau} M(s,\tau) A e^{-(\tau-s)A} A \hat{k}_{n}(s) ds = P_1(s)+ P_2(s),
\end{align*}
where
\begin{align*}
		P_1(s) = \int_0^{\frac{\tau}{2}} \frac{s^2}{2} A e^{-(\tau-s)A}A \hat{k}_{n}(s) ds, \qquad P_2(s) = \int_{\frac{\tau}{2}}^{\tau} \frac{(\tau-s)^2}{2} A e^{-(\tau-s)A}A \hat{k}_{n}(s) ds.
\end{align*}
Since $A \hat{k}_{n}(0)= \mathcal{O}(1)$, by using \eqref{B3.2} and the parabolic smoothing property \eqref{B3.3}, we get
\begin{equation}\label{B3.30}
	\begin{aligned}
		P_1(s) &=  \int_0^{\frac{\tau}{2}}  \frac{s^2}{2}A e^{-(\tau-s)A}A \hat{k}_{n}(s) ds  \\
		& = \int_0^{\frac{\tau}{2}}  \frac{s^2}{2} A e^{-(\tau-s)A}A \hat{k}_{n}(0) ds + \int_0^{\frac{\tau}{2}}\frac{s^2}{2}  A^2 e^{-(\tau-s)A} \mathcal{O}(s) ds\\
		& =A \int_0^{\frac{\tau}{2}} \frac{s^2}{2} e^{-(\tau-s)A} \left(A \hat{k}_{n}(0) \right) ds  \\
		& \qquad + A \left( A e^{-\frac{\tau}{2}A} \right) \int_0^{\frac{\tau}{2}} \frac{s^2}{2}
		e^{-\left( \frac{\tau}{2}-s\right) A} \mathcal{O}(s) ds\\
		& = A \mathcal{O}(\tau^3),
	\end{aligned}
\end{equation}
and
\begin{equation}\label{B3.31}
	\begin{aligned}
		P_2(s) &=  \int_{\frac{\tau}{2}}^{\tau} \frac{(\tau-s)^2}{2} A e^{-(\tau-s)A}A \hat{k}_{n}(s) ds \\
		& = A\int_{\frac{\tau}{2}}^{\tau} \frac{(\tau-s)^2}{2}  e^{-(\tau-s)A}  \left(A \hat{k}_{n}(0) \right)ds\\
		& \qquad + A \int_{\frac{\tau}{2}}^{\tau}(\tau-s) \Big( (\tau-s) e^{-(\tau-s)A} A \Big) \mathcal{O}(s) ds\\
		& = A \mathcal{O}(\tau^3).
	\end{aligned}
\end{equation}
Further, by \eqref{B3.2} and Assumption \ref{Bas3}, $\hat{k}_{n}$ is twice continuously differentiable. Therefore, we get the following estimates:
\begin{equation}\label{B3.33}
	\begin{aligned}
		Q(s) =  2 \int_0^{\tau} M(s,\tau) Ae^{-(\tau-s)A} \hat{k}_{n}^{\prime}(s) ds = A \mathcal{O}(\tau^3),
	\end{aligned}
\end{equation}
and
\begin{equation}\label{B3.34}
	\begin{aligned}
		R(s) =   \int_0^{\tau} M(s,\tau) e^{-(\tau-s)A} \hat{k}_{n}^{\prime \prime}(s) ds  = \mathcal{O}(\tau^3).
	\end{aligned}
\end{equation}
Combining \eqref{B3.30}, \eqref{B3.31}, \eqref{B3.33} and \eqref{B3.34}, we obtain the following representation of $\delta_{n+1}^{\left[1,2\right] }$ :
\begin{equation}
	\begin{aligned}
		\delta_{n+1}^{\left[1,2\right] } = A \mathcal{O}(\tau^3) + \mathcal{O}(\tau^3).
	\end{aligned}
\end{equation}
This together with \eqref{B3.28} yields that
\begin{align}\label{B3.36}
	\delta_{n+1}^{\left[ 1\right] } =  	\delta_{n+1}^{\left[1,1\right]}+ \delta_{n+1}^{\left[1,2\right]} =  A\mathcal{O}(\tau^3) + \mathcal{O}(\tau^3).
\end{align}


To complete the local error estimate, we still need to bound $\| \delta_{n+1}^{\left[2 \right] }\|$ in \eqref{B3.17}. We set
\begin{align*}
	l_{n}(s) = e^{-(\tau-s)A} \hat{l}_{n}(s) ,\quad \hat{l}_{n}(s)= h(t_{n}+s, \hat{v}(t_{n}+s)) - h(t_{n}+s, \hat{u}(t_{n}+s)),
\end{align*}
and
\begin{align*}
	\tilde{l}_{n}(s)=  e^{-(\tau-s)A}  \bar{l}_{n}(s),\quad   \bar{l}_{n}(s)  = h(t_{n}+s, \bar{v}(t_{n}+s)) - h(t_{n}+s, \hat{u}(t_{n}+s)).
\end{align*}
Thus, we can rewrite $ \delta_{n+1}^{\left[2 \right] }$ as follows:
\begin{align*}
	\delta_{n+1}^{\left[2 \right] } = \int_0^{\frac{\tau}{2}} {l}_{n}(s) ds+ \int_{\frac{\tau}{2}}^{\tau} \tilde{l}_{n}(s) ds.
\end{align*}
Taylor expansion of $l_{n}$ leads to
\begin{align*}
	l_{n}(s)= l_{n}\left( \frac{\tau}{2} \right) + \mathcal{O}(\tau^2).
\end{align*}
In view of \eqref{B3.22}, \eqref{B3.24}, \eqref{B3.26}, and by using $ \nabla \hat{u}(t_{n}) = 0$, we get
\begin{align}\label{B3.37}
	l_{n}\left( \frac{\tau}{2} \right) = \mathcal{O}(\tau^2), \qquad	\int_0^{\frac{\tau}{2}} {l}_{n}(s) ds	= \mathcal{O}(\tau^3).
\end{align}
In a similar way, using \eqref{B3.11}, \eqref{B3.24}, \eqref{B3.26}, and $\nabla \hat{u}(t_{n}) = 0$, we obtain
\begin{align}\label{B3.38}
\tilde{l}_{n}(s) = \mathcal{O}(\tau^2), \qquad	\int_{\frac{\tau}{2}}^{\tau} \tilde{l}_{n}(s) ds	= \mathcal{O}(\tau^3).
\end{align}
Combining \eqref{B3.37} and \eqref{B3.38} we arrive at
\begin{align}\label{B3.39}
	\delta_{n+1}^{\left[2 \right] } = \mathcal{O}(\tau^3).
\end{align}
From \eqref{B3.15}, combining \eqref{B3.36} and \eqref{B3.39}, we get the following representation of the local error
\begin{align}\label{B3.42}
	\delta_{n+1} = 	\delta_{n+1}^{\left[1 \right] }+ \delta_{n+1}^{\left[2 \right] }=
	A \mathcal{O}(\tau^3)+ \mathcal{O}(\tau^3).
\end{align}

Next, we study the global error of the modified Strang splitting method applied to~\eqref{B1.1}.

\subsection{Global error}
In this subsection, we show that the modified Strang splitting method is convergent of order two. In particular, we provide the proof of Theorem $\ref{Th: 2.1}$. Let us denote by $e_{n}=u_{n}-u(t_{n})$ the global error. Thus, we have
\begin{align}\label{B3.43}
	e_{n+1}=\mathcal{S}_{\tau}u_{n}-\mathcal{S}_{\tau} u(t_{n})+\delta_{n+1},
\end{align}	
where $\delta_{n+1}= \mathcal{S}_{\tau} u(t_{n})- u(t_{n+1})$ denotes the local error. The numerical solution $S_\tau u_n$ can be expressed in the same way as \eqref{B3.14} by first integrating \eqref{B3.7} with the initial value $\hat{v}(t_n)=0$ and then proceeding as described in Algorithm~\ref{alg 1} or~\ref{alg 2}. In order to distinguish the arising functions from those in \eqref{B3.14}, we call them $\hat{v}_\text{nu}(t)$, $\bar{v}_\text{nu}(t)$, and
$$
Y_\text{nu} = \int_0^{\frac{\tau}{2}} e^{ -\left( \frac{\tau}{2}-s \right)A } h(t_{n}+s, \hat{v}_\text{nu}(t_{n}+s)) ds.
$$
This shows that
\begin{equation}\label{B3.44}
	\begin{aligned}
		S_{\tau}u_{n}& = z(t_{n+1}) + \int_0^{\frac{\tau}{2}} e^{-(\tau-s)A} h(t_{n}+s, \hat{v}_\text{nu}(t_{n}+s)) ds \\
		& \qquad + \tau e^{-\frac{\tau}{2} A} a \cdot \nabla Y_\text{nu} + \frac{\tau^2}{2} e^{-\frac{\tau}{2} A} a \cdot \nabla \left( a \cdot \nabla Y_\text{nu} \right)  \\
		& \qquad + \int_{\frac{\tau}{2}}^{\tau} e^{-(\tau-s)A} h(t_{n}+s, \bar{v}_\text{nu}(t_{n}+s)) ds + \mathcal{O}(\tau^3).
	\end{aligned}
\end{equation}
Taking the difference of \eqref{B3.44} and \eqref{B3.14} results in
\begin{align}\label{B3.46}
	e_{n+1}= e^{-\tau A} e_{n} + E_n + \delta_{n+1}+\mathcal{O}(\tau^3),
\end{align}
where $E_n$ are the differences of the corresponding terms in \eqref{B3.44} and \eqref{B3.14}, respectively. Solving \eqref{B3.46} and using $\| e_0\| =0$ together with the estimate \eqref{B3.42} of the defects and the Lipschitz continuity of $f$ gives
\begin{equation}\label{B3.47}
\begin{aligned}
	\|e_{n}\| & \le  C \tau  \sum_{k=1}^{n-1} \left\|  e^{-(n-k-1)\tau A}\right\| \left\| e_{k}\right\|+C \tau \sum_{k=1}^{n-1} \left\|  e^{-\left(n-k-\frac12\right)\tau A} a \cdot \nabla \right\| \left\|e_{k}\right\| \\
	& \qquad +	\sum_{k=0}^{n-2} \left\| e^{-(n-k-1)\tau A} \left(A \mathcal{O}(\tau^3)+ \mathcal{O}(\tau^3)\right) \right\|  + \| \delta_{n} \| + \mathcal{O}(\tau^2).
\end{aligned}
\end{equation}
Using the parabolic smoothing \eqref{B3.3} and the fact that $\| \delta_{n}\| = \mathcal{O}(\tau^2)$, we obtain
\begin{equation}
	\begin{aligned}
		\| e_{n}\| & \le C \tau \sum_{k=1}^{n-1} \| e_{k}\| + C	\tau \sum_{k=1}^{n-1} \left\| e^{-(n-k)\tau A} A^{\frac{1}{2}} \right\| \left\| A^{-\frac{1}{2}} a \cdot\nabla \right\|  \| e_{k}\| \\
		& \qquad + 	C\tau^3 \sum_{k=0}^{n-2} \left\| e^{-(n-k-1)\tau A} A\right\| + nC \tau^3  + C\tau^2 \\
		&\le    C \tau \sum_{k=1}^{n-1} \| e_{k}\| + C\tau \sum_{k=1}^{n-1}t_{n-k}^{-\frac{1}{2}} \Vert e_{k}\Vert + C\tau^3 \sum_{k=1}^{n-1} \frac{1}{k\tau} + C\tau^2 \\
		& \le C \tau \sum_{k=1}^{n-1} \| e_{k}\| + C\tau \sum_{k=1}^{n-1}  t_{n-k}^{-\frac{1}{2}} \Vert e_{k}\Vert+ C\tau^2(1+ \log \left| \tau\right| ).
	\end{aligned}
\end{equation}
The application of a discrete Gronwall lemma completes the proof of the second-order convergence of the modified Strang splitting.
\section{Numerical results in one space dimension}
In this section, we provide a series of numerical results for the convection-diffusion-reaction problem \eqref{B1.1} on the domain $\Omega=[0,1]$ and with $D=\frac1{10} \partial_{xx}$ and $a(x)=x^2$. The boundary conditions on the left and right sides are denoted by $b_1$ and $b_2$, respectively. We employ the standard centered second-order finite differences to discretize the Laplacian and the first-order upwind scheme for the convection term. For the discretization we use 200 grid points.
The reference solution is computed using the ODE45 solver with absolute and relative tolerances set to $10^{-9}$. In the simulations presented here, the scheme \eqref{B1.4} is referred to as the classical Strang splitting, while the scheme given by Algorithm \ref{alg 1} or \ref{alg 2} (an initial-corrected splitting method) is referred to as the modified Strang splitting.

\begin{figure}[ht]
		\centering
		\begin{subfigure}[b]{0.49\textwidth}
			\centering
			\includegraphics[width=\textwidth]{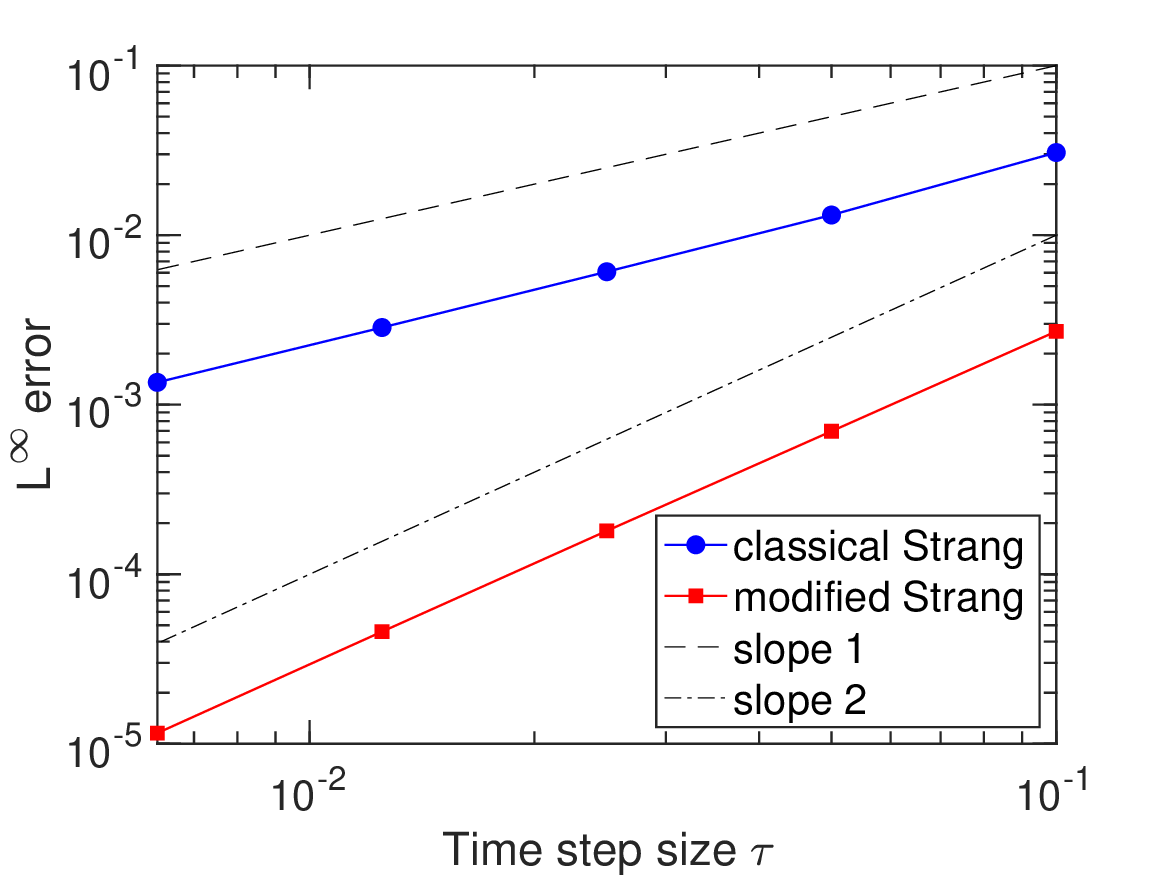}
			\caption{Homogeneous BCs}
			\label{fig:1a}
		\end{subfigure}
		\hfill
		\begin{subfigure}[b]{0.49\textwidth}
			\centering
			\includegraphics[width=\textwidth]{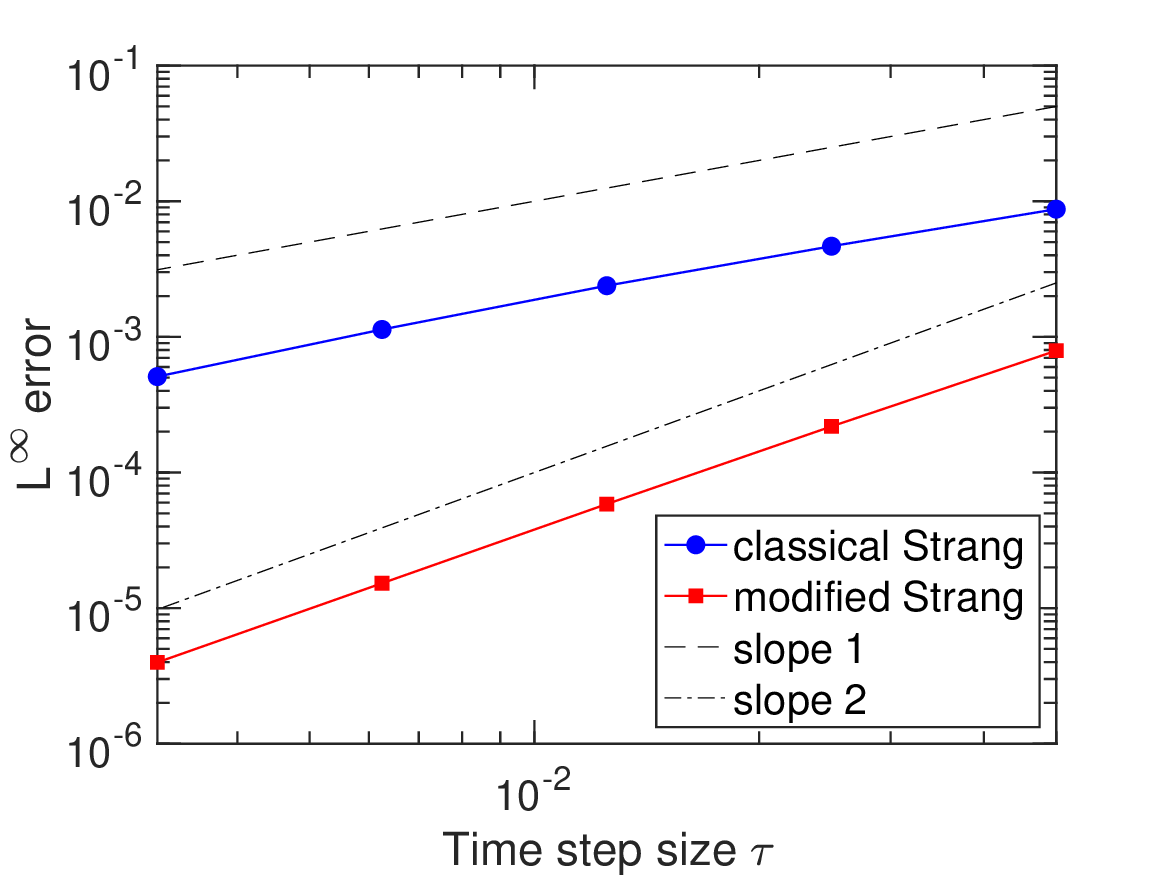}
			\caption{Inhomogeneous BCs}
			\label{fig:1b}
		\end{subfigure}
	\hfill
	\begin{subfigure}[b]{0.49\textwidth}
		\centering
		\includegraphics[width=\textwidth]{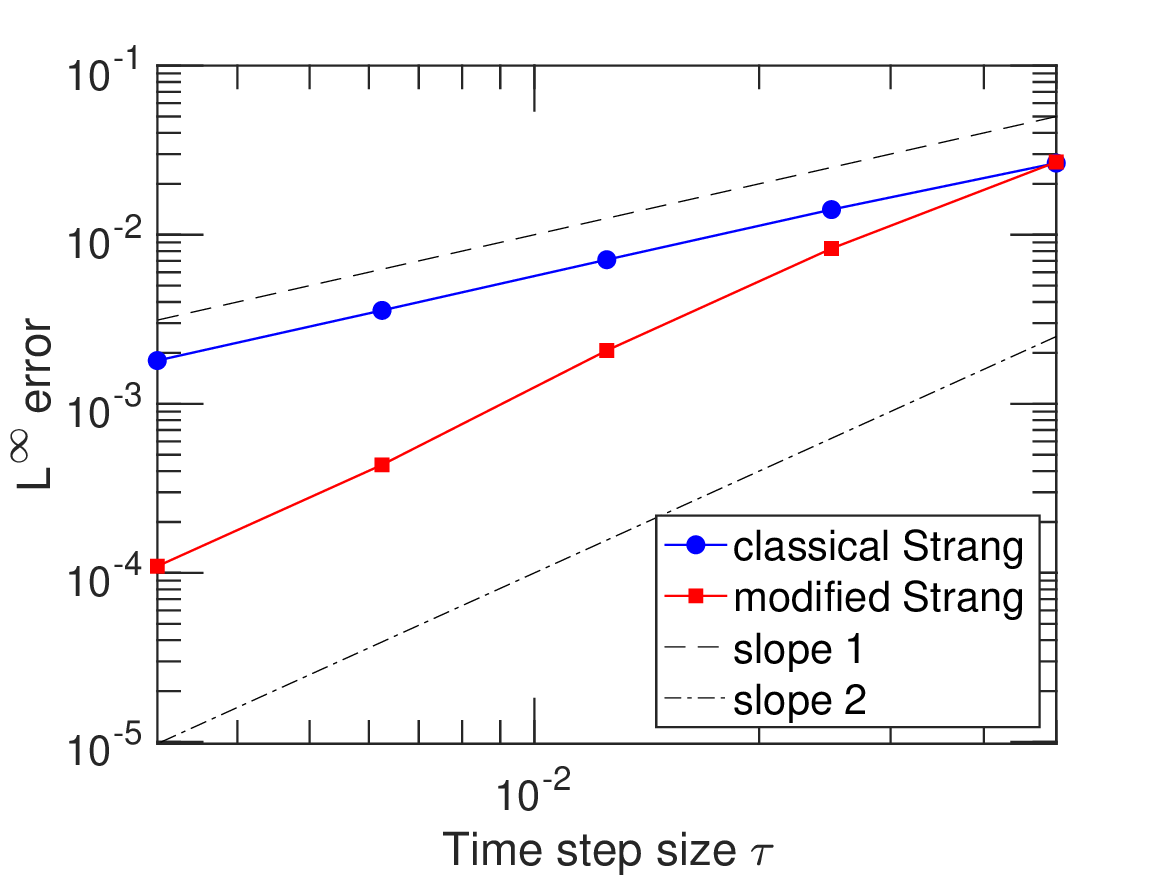}
		\caption{Time-dependent BCs}
		\label{fig:1c}
	\end{subfigure}
		\caption{The absolute error in the discrete infinity norm is computed at $t=1$ by comparing the numerical solution to a reference solution.}
		\label{fig1}
\end{figure}

\begin{exmp}
	In this example, we consider \eqref{B1.1} with the reaction term $f(t, u(t,x))= u^2(t,x)+ \phi(t,x)$ where the source function $\phi(t,x)$ is chosen such that $u(t,x)= x(1-x)e^{t} $ is the exact solution. The source term $\phi(t,x)$ is given by
	\begin{align*}
		\phi(t,x)= (1-3x-x^2)e^{t}- (x^2 -2x^3+x^4)e^{2t}.
	\end{align*}
	We impose homogeneous Dirichlet boundary conditions, i.e., $b_1 =b_2 =0$. The initial data $u_0(x)= x(1-x)$ is chosen to comply with these boundary conditions. The numerical results are displayed in Figure \ref{fig:1a}.
\end{exmp}
\begin{exmp}
This example considers \eqref{B1.1} with the reaction term $f(t, u(t,x)) = u^2(t,x)$, subject to inhomogeneous Dirichlet boundary conditions with $b_1 = 1$ and $b_2=2$. The initial data is chosen as $u_0(x)= 1 + \sin\left( \frac{\pi}{2}x\right)$, ensuring that the prescribed boundary conditions are satisfied. The numerical results are shown in Figure \ref{fig:1b}.
\end{exmp}

\begin{exmp}
In this example, we consider \eqref{B1.1} with the reaction term $f(t,x)$ defined
	\begin{align*}
		f(t,x) = e^{t} \Big( 1 + (1+ 0.2 \pi^2) \sin^2(\pi x) - 0.2 \pi^2 \cos^2(\pi x) - x^2 \pi \sin(2 \pi x)\Big).
	\end{align*}
The problem is subject to time-dependent boundary conditions, with the left boundary condition $b_1(t)=1 + \sin(5t) $ and the right boundary condition $b_2(t) = 1+ \sin(10 \pi t)$. The initial data is given by $u_0(x) = 1 + \sin^2(\pi x)$. The numerical results are shown in Figure \ref{fig:1c}.	
\end{exmp}

It is observed that all simulations presented in Figure \ref{fig1} show the expected order reduction to approximately order one for the classical Strang splitting scheme, even when homogeneous Dirichlet boundary conditions are imposed. In contrast, the modified Strang splitting method achieves second-order convergence, exhibiting significantly enhanced accuracy regardless of whether homogeneous, inhomogeneous or time-dependent Dirichlet boundary conditions are applied.		

\section{Numerical results in two space dimensions}
\begin{figure}[ht]
	\centering
	\includegraphics[width=0.5\linewidth]{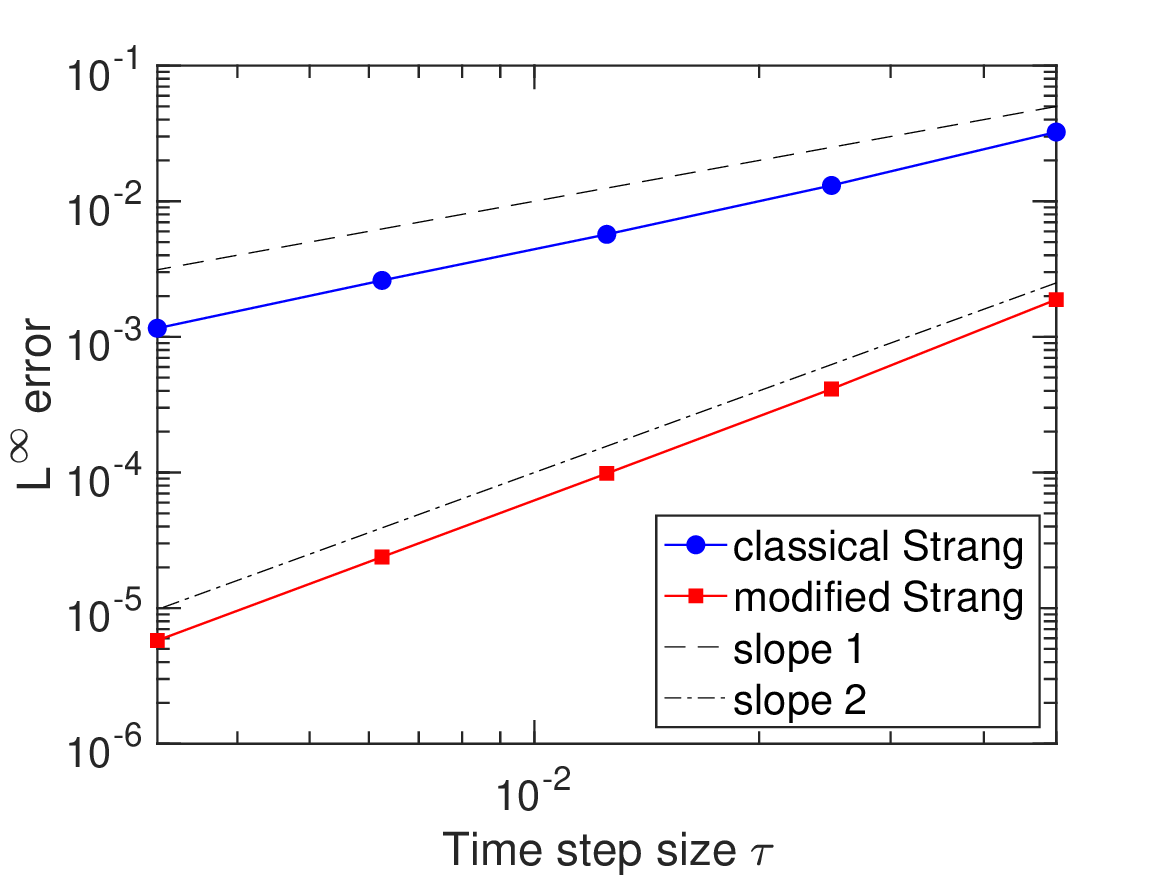}
	\caption{The absolute error in the discrete infinity norm is computed at $t=1$ by comparing the numerical solution to a reference solution. Slope $2 $ is displayed by a dash-dotted line.}
	\label{fig3}
\end{figure}

In this section, we consider the convection-diffusion-reaction problem \eqref{B1.1} with $f(t, u(t)) = e^{u(t)}$ on $\Omega = [0,1]^2$, where $D$ is the standard second-order finite difference approximation of the Laplacian, using $50 \times 50$ grid points. The diffusion coefficient is set to $0.1$, while the velocity field is defined as $a = (x^2, y^2)$. Homogeneous Dirichlet boundary conditions are imposed. The initial condition is specified as
\begin{align*}
	u_0(x,y) = e^{-100(x-0.5)^2} e^{-100(y-0.5)^2} \sin(\pi x)\sin(\pi y).
\end{align*}
The reference solution is computed using ODE45 solver, with the absolute and relative tolerances are set to $10^{-7}$.
The error in the discrete infinity norm is computed at $t= 1$ by comparing the numerical solution to a reference solution.
The numerical results are displayed in Figure \ref{fig3}.
We observe an order reduction to approximately order one in the infinity norm for the classical Strang splitting, whereas for the modified Strang splitting, second-order accuracy is observed. 	

\section{Conclusion}
We have proposed an initial-corrected splitting method that effectively overcomes the observed order reduction when splitting convection-diffusion-reaction problems with Dirichlet boundary conditions. Specifically, we have proposed a correction for which we have demonstrated second-order convergence both analytically and numerically. Moreover, the proposed approach is easy to apply and is computationally efficient.

\section*{Acknowledgements}
This project has received funding from the European Union’s Horizon 2020 research and innovation programme under the Marie Sk\l{}odowska-Curie grant agreement No 847476.

\bibliography{p2_references}
\nocite{*}

\end{document}